\documentclass[11pt,a4paper]{article}

\usepackage[T1]{fontenc}
\usepackage{lmodern}
\usepackage{amsmath,amssymb,amsthm}
\usepackage[margin=2.6cm]{geometry}
\usepackage[hidelinks]{hyperref}
\usepackage{authblk}
\usepackage{float}
\usepackage{booktabs}

\theoremstyle{plain}
\newtheorem{theorem}{Theorem}
\newtheorem{example}{Example}

\title{Lean-verified lower bounds for the Shannon capacity of odd~cycles}
\author[1]{Pjotr Buys}
\author[2]{Sven Polak}
\author[3]{Jeroen Zuiddam}
\date{}
\affil[1]{Centrum Wiskunde \& Informatica}
\affil[2]{Tilburg University}
\affil[3]{University of Amsterdam}

\begin{document}
\maketitle

\begin{abstract}
\begin{sloppypar}
We give new lower bounds for the Shannon capacities of small odd cycles:
$\Theta(C_7)\geq3.258805369885\ldots$,
$\Theta(C_{11})\geq5.294502522149\ldots$,
$\Theta(C_{13})\geq6.302455083464\ldots$,
$\Theta(C_{15})\geq7.301600534487\ldots$,
$\Theta(C_{19})\geq9.357192705918\ldots$,
$\Theta(C_{21})\geq10.342455853338\ldots$, and
$\Theta(C_{23})\geq11.328224257774\ldots$.
The bounds are obtained by an iterative procedure due to
Gao (2026) which is based on a method by Itty, Rosin, Carstensen and Reichman (2026). The bounds are fully formalised in Lean.
\end{sloppypar}
\end{abstract}

\section{Introduction}

For a finite graph~$G$, let~$G^{\boxtimes d}$ denote its $d$th strong power and let~$\alpha(G^{\boxtimes d})$ denote its independence number. The \emph{Shannon capacity} of~$G$ is
$
\Theta(G):=\sup_{d\geq1}\alpha(G^{\boxtimes d})^{1/d}.
$
Shannon capacity was introduced by Shannon in zero-error information
theory~\cite{shannon}. Lov\'asz introduced the theta function $\vartheta(G)$ as an upper bound on Shannon capacity $\Theta(G)$ and used it to determine $\Theta(C_5)=\sqrt5$~\cite{lovasz}. The
capacity of every odd cycle of length at least seven has remained open since.
Classical combinatorial packing constructions for powers of odd cycles go back
to Baumert, McEliece, Rodemich, Rumsey, Stanley and Taylor~\cite{baumert}.

Recently, \cite{itty, gao} used large language models to obtain improved lower bounds on Shannon capacity of small odd cycles. Following their methods and  using ChatGPT 5.6 Sol Pro and Claude Opus 5 we find further improved bounds.

\section{Valid tuples and the product theorem}

We use the following reformulation of Gao's gadget construction
and product lemma~\cite[Definitions~2--3 and Lemma~4]{gao}.
Let \(G\) be a finite graph. A \emph{valid tuple} for \(G\) is
\(\tau=(I,S,f_0,f_1,X)\), where \(I,S,X\subseteq V(G)\) and
\(f_0,f_1:S\to V(G)\), such that \(I\) is independent, \(S\subseteq I\),
and, for every \(s\in S\), exactly one of \(f_0(s)\) and \(f_1(s)\)
equals \(s\), while \(N_G[f_i(s)]\cap I=\{s\}\) for \(i\in\{0,1\}\).
Furthermore, \(f_0(S)\), \(f_1(S)\), and \(X\) are independent, and
\(X\cap N_G[f_0(S)]\cap N_G[f_1(S)]=\varnothing\).
The \emph{profile} of \(\tau\) for $G$ is
\[
    \Pi(\tau)=
    \left(
        |I|,\ |S|,\ |X|,\
        \left|X\setminus
        \bigl(N_G[f_0(S)]\cup N_G[f_1(S)]\bigr)\right|
    \right).
\]
For two quadruples of integers, write
\((a_1,b_1,c_1,d_1)\star(a_2,b_2,c_2,d_2)\) for
\[
    \left(
        (a_1-b_1)(a_2-b_2)+c_1b_2+b_1c_2,\ 
        d_1b_2+b_1d_2,\ 
        c_1c_2,\ 
        d_1d_2+(c_1-d_1)(c_2-d_2)
    \right).
\]

\begin{theorem}[\cite{gao}]
Let \(G\) and \(H\) be finite graphs. If \(\tau\) is a valid tuple for
\(G\) and \(\tau'\) is a valid tuple for \(H\), then there exists a valid
tuple \(\tau\star\tau'\) for \(G\boxtimes H\) such that
\(\Pi(\tau\star\tau')=\Pi(\tau)\star\Pi(\tau')\).
\end{theorem}

The first coordinate of \(\Pi(\tau)\) is \(|I|\), so every valid tuple
directly supplies an independent set of that size. The product operation
can produce independent sets in strong powers that are larger than the
ordinary products of the original independent sets. Iterating it can
therefore improve lower bounds on Shannon capacity.

\begin{example}
On \(C_5=\mathbb Z_5\), take \(I=\{0,2\}\), \(S=\{0\}\),
\(f_0(0)=0\), \(f_1(0)=4\), and \(X=\{1,3\}\). This gives a valid tuple
\(\tau\) with \(\Pi(\tau)=(2,1,2,0)\). Hence
\(\Pi(\tau\star\tau)=(5,0,4,4)\), and therefore
\(\alpha(C_5^{\boxtimes 2})\geq 5\).
\end{example}

\section{Bounds obtained by iteration}

In this section we describe new lower bounds on the Shannon capacities
of several small odd cycles. In each case we start from a
computationally obtained valid tuple in a fixed power. We record only
its profile and the sequence of \(\star\)-products leading to the
stated bound. The base valid tuples
are recorded in the accompanying Lean formalisation, where they
are also verified.

\begin{table}[H]
\centering
\begin{tabular}{@{}rrllll@{}}
\toprule
$n$ & power & new bound & previous bound & improvement & $\vartheta(C_n)$ \\
\midrule
$7$  & $200$   & $3.2588053$  & $3.2587891$~\cite{gao}      & $1.6\cdot10^{-5}$ & $3.3176672$ \\
$11$ & $174$   & $5.2945025$  & $5.2897736$~\cite{itty}     & $4.7\cdot10^{-3}$ & $5.3863029$ \\
$13$ & $432$   & $6.3024550$  & $6.3001091$~\cite{itty}     & $2.3\cdot10^{-3}$ & $6.4041685$ \\
$15$ & $4096$  & $7.3016005$  & $7.3013990$~\cite{itty}     & $2.0\cdot10^{-4}$ & $7.4171482$ \\
$19$ & $16384$ & $9.3571927$  & $9.3571200$~\cite{baumert}  & $7.3\cdot10^{-5}$ & $9.4347713$ \\
$21$ & $4096$  & $10.3424558$ & $10.3422729$~\cite{baumert} & $1.8\cdot10^{-4}$ & $10.4410325$ \\
$23$ & $2048$  & $11.3282242$ & $11.3278379$~\cite{baumert} & $3.9\cdot10^{-4}$ & $11.4461936$ \\
\bottomrule
\end{tabular}
\caption{Lower bounds on $\Theta(C_n)$ obtained here, the strong power in which
the independent set witnessing the bound lives, the previously best known bound,
and the Lov\'asz bound $\vartheta(C_n)=n\cos(\pi/n)/(1+\cos(\pi/n))$. Values are
truncated, not rounded.}
\label{tab:bounds}
\end{table}

\subsection*{The cycle \(C_7\)}

We start with a valid tuple $\tau$ for \(C_7^{\boxtimes 5}\) whose profile $\Pi(\tau)$ is
\(\pi_1=(367,8,367,322)\). Form
\(\pi_2=\pi_1\star\pi_1\), \(\pi_3=\pi_1\star\pi_2\),
\(\pi_5=\pi_2\star\pi_3\), and then
\(\pi_{2r}=\pi_r\star\pi_r\) for \(r=5,10,20\). The resulting \(\pi_{40}\)
is a profile for \(C_7^{\boxtimes 200}\) and gives
\(\Theta(C_7)\geq 3.258805369885\ldots\).

Lower bounds for $C_7$ were obtained in~\cite{baumert,veszer,matost}. A size $367$ independent set in $C_7^{\boxtimes 5}$ which improved the lower bound was discovered by Polak and Schrijver~\cite{polsch}, which was later recovered by
FunSearch~\cite{funsearch}.  Recently, Itty,
Rosin, Carstensen and Reichman~\cite{itty} improved the lower bound using a product construction in~$C_7^{\boxtimes 10}$, which was subsequently improved by Gao using a recursive construction~\cite{gao} to $\Theta(C_7) \geq 3.258789153908$.

\subsection*{The cycle \(C_{11}\)}

We start with a valid tuple for \(C_{11}^{\boxtimes 3}\) whose profile is
\(\pi_1=(148,3,148,142)\). Form
\(\pi_2=\pi_1\star\pi_1\), \(\pi_3=\pi_1\star\pi_2\),
\(\pi_4=\pi_2\star\pi_2\), \(\pi_7=\pi_3\star\pi_4\),
\(\pi_8=\pi_4\star\pi_4\), \(\pi_{14}=\pi_7\star\pi_7\),
\(\pi_{15}=\pi_7\star\pi_8\),
\(\pi_{29}=\pi_{14}\star\pi_{15}\), and
\(\pi_{58}=\pi_{29}\star\pi_{29}\). The last profile lies in
\(C_{11}^{\boxtimes 174}\) and gives
\(\Theta(C_{11})\geq 5.294502522149\ldots\).

The classical independent set of size~$148$ in $C_{11}^{\boxtimes 3}$ is due to Baumert et
al.~\cite{baumert}; Itty et al.\ subsequently improved the lower bound on
$\Theta(C_{11})$~\cite{itty} to $\Theta(C_{11}) \geq 5.289773694291$.

\subsection*{The cycle \(C_{13}\)}

We start with a valid tuple for \(C_{13}^{\boxtimes 6}\) whose profile is
\(\pi_1=(62530,1014,62530,60502)\). Form
\(\pi_2=\pi_1\star\pi_1\), \(\pi_3=\pi_1\star\pi_2\),
\(\pi_4=\pi_2\star\pi_2\), \(\pi_5=\pi_2\star\pi_3\),
\(\pi_9=\pi_4\star\pi_5\), and then
\(\pi_{2r}=\pi_r\star\pi_r\) for \(r=9,18,36\). The resulting \(\pi_{72}\)
lies in \(C_{13}^{\boxtimes 432}\) and gives
\(\Theta(C_{13})\geq 6.302455083464\ldots\).

The size $247$ cube construction goes back to Baumert et
al.~\cite{baumert}, and its optimality was established by Bohman,
Holzman and Natarajan~\cite{bohman}; Itty et al.\ later improved the
capacity lower bound~\cite{itty} to $\Theta(C_{13}) \geq 6.300109117130$.

\subsection*{The cycle \(C_{15}\)}

We start with a valid tuple for \(C_{15}^{\boxtimes 4}\) whose profile is
\(\pi_1=(2842,3,2842,2833)\). Form
\(\pi_{2r}=\pi_r\star\pi_r\) successively for
\(r=1,2,4,\ldots,512\). The resulting \(\pi_{1024}\) lies in
\(C_{15}^{\boxtimes 4096}\) and gives
\(\Theta(C_{15})\geq 7.301600534487\ldots\).

A lower bound for $C_{15}$ was obtained by Mathew and
\"Osterg{\aa}rd~\cite{matost}. Thereafter, De~Boer, Buys and Zuiddam~\cite{deboer} proved that~$\Theta(C_{15}) \geq \sqrt[4]{2842}$, after which Itty et al.~\cite{itty} proved that $\Theta(C_{15}) \geq 8076974^{1/8} = 7.301399060577\ldots$.

\subsection*{The cycle \(C_{19}\)}

We start with a valid tuple for \(C_{19}^{\boxtimes 4}\) whose profile is
\(\pi_1=(7666,2,7666,7661)\). Form
\(\pi_{2r}=\pi_r\star\pi_r\) successively for
\(r=1,2,4,\ldots,2048\). The resulting \(\pi_{4096}\) lies in
\(C_{19}^{\boxtimes 16384}\) and gives
\(\Theta(C_{19})\geq 9.357192705918\ldots\).

\subsection*{The cycle \(C_{21}\)}

We start with a valid tuple for \(C_{21}^{\boxtimes 4}\) whose profile is
\(\pi_1=(11441,10,11441,11398)\). Form
\(\pi_{2r}=\pi_r\star\pi_r\) successively for
\(r=1,2,4,\ldots,512\). The resulting \(\pi_{1024}\) lies in
\(C_{21}^{\boxtimes 4096}\) and gives
\(\Theta(C_{21})\geq 10.342455853338\ldots\).

\subsection*{The cycle \(C_{23}\)}

We start with a valid tuple for \(C_{23}^{\boxtimes 4}\) whose profile is
\(\pi_1=(16466,30,16466,16323)\). Form
\(\pi_{2r}=\pi_r\star\pi_r\) successively for
\(r=1,2,4,\ldots,256\). The resulting \(\pi_{512}\) lies in
\(C_{23}^{\boxtimes 2048}\) and gives
\(\Theta(C_{23})\geq 11.328224257774\ldots\).

\section{Lean formalisation}

All the Shannon capacity bounds stated above are formalised in the proof assistant Lean and can be found in the repository: \begin{center}
\url{https://github.com/spectra-research/shannon-capacity-lean}
\end{center}
For each of the seven base tuples, the finite sets $I,S,X$ and the maps
$f_0,f_1$ appear as explicit literals in the Lean sources
\texttt{ShannonBounds/BaseC}$n$\texttt{Data.lean}, together with the proofs
that each tuple is valid. These are the same literals that the Lean
formalisation checks.

\bibliographystyle{alphaurl}
\bibliography{refs}

\end{document}